\documentclass[11pt]{article}
\usepackage{lineno}
\modulolinenumbers[1]
\usepackage{amssymb,graphicx,amssymb,amsmath,amsopn, mathtools,color}
\definecolor{brightmaroon}{rgb}{0.76, 0.13, 0.28}
\definecolor{airforceblue}{rgb}{0, 0.25, 0.77}
\usepackage[bookmarks=true,colorlinks=true,linkcolor=airforceblue,citecolor=brightmaroon,plainpages=tfalse,pdftex]{hyperref}
\hypersetup{bookmarksopen=true,pdfview=fitbh}

\usepackage{cases}
\usepackage{textcomp}
\usepackage{tikz-cd}
\usepackage{amsthm}
\usepackage[mathscr]{euscript}
\usepackage{amscd}
\usepackage{pb-diagram}

\allowdisplaybreaks
\newtheorem{pro}{Proposition}
\newtheorem{defi}[pro]{Definition}
\newtheorem{teo}[pro]{Theorem}

\newtheorem{coro}[pro]{Corollary}

\newtheorem{rem}[pro]{Remark}
\swapnumbers
\newtheorem{exa}[]{Example}

\newcommand{\N}{\mathbb{N}}

\newcommand{\Su}{\boldsymbol{\mathscr{S}}}
\newcommand{\un}{\mathbf{u}}
\newcommand{\vn}{\mathbf{v}}

\newcommand{\prodint}[1]{\left\langle{#1}\right\rangle}
\newcommand*{\Scale}[2][4]{\scalebox{#1}{$#2$}}
\newcommand{\deltan}{\boldsymbol{\delta}}

\title{Associated  orthogonal polynomials of the first kind  and  Darboux  transformations }

\author{J. C. Garc\'ia-Ardila$^1$, F. Marcell\'an$^{2}$, P. H. Villamil-Hern\'andez $^{2}$ \\
	$^1$ Departamento de Matem\'atica Aplicada a la Ingenier\'ia Industrial, \\Universidad Polit\'ecnica de Madrid,\\ Calle Jos\'e Gutierrez Abascal 2, 28006 Madrid, Spain.\\
	$^2$Departamento de Matem\'aticas, Universidad Carlos III de Madrid,\\ Avenida  Universidad 30, 28911 Leganés, Spain\\}

\begin{document}

\maketitle
 	\begin{abstract}
 	
 	Let $\un,$ be a quasi-definite linear functional defined on the space of polynomials $\mathbb{P}.$	For such a functional we can define  a sequence  of monic orthogonal polynomials (SMOP in short) $(P_n)_{n\geq 0},$ which satisfies a three term recurrence relation. Shifting one unity the recurrence coefficient indices given the sequence of associated polynomials of the first kind which are orthogonal with respect to a linear functional denoted by $\un^{(1)}$.
 	
 	In the literature two special transformations of the functional $\un$ are studied, the canonical Christoffel transformation   $\widetilde \un=(x-c)\un$ and the canonical Geronimus transformation $\widehat \un=\frac{\un}{(x-c)}+M\deltan_c$ ,  where $c$ is a fixed  complex number,  $M $ is  a free parameter  and $\deltan_c$ is the linear  functional defined on $\mathbb{P}$  as $\prodint{\deltan_c,p(x)}=p(c).$
 	
 	For the Christoffel transformation with SMOP $(\widetilde P_n)_{n\geq 0}$, we are interested in analyzing  the relation between the linear functionals  $ \un^{(1)}$ and $\widetilde{\un}^{(1)}.$ There,  the super index denotes the linear functionals associated with the  orthogonal polynomial sequences of the first kind $(P_n^{(1)})_{n\geq 0}$ and $(\widetilde P_n^{(1)})_{n\geq 0},$ respectively. This problem is also studied for Geronimus transformations. Here we give close relations between their corresponding monic Jacobi matrices by using the LU and UL factorizations. To get this result, we first need to study the relation between $\un^{-1}$ (the inverse functional) and $\un^{(1)}$ which can be expressed  from a quadratic Geronimus transformation.
 \end{abstract}


\section{Introduction and basic background}

First of all we will introduce the basic background about linear functionals that will be used in the sequel in order to have a self-contained presentation. The reader can check \cite{Maro1991}
where  an overview containing a more detailed analysis of algebraic properties of linear functionals, which are the main tools to deal with the standard theory of orthogonal polynomials, is presented.

Let $\un$ be a complex-valued linear functional   defined on the linear space of polynomials with complex coefficients $\mathbb{P}$, i.e. $\un:\mathbb{P}\to \mathbb{C}$, $p(x)\to \prodint{\un,p(x)}$.  We denote the $n$-th moment of $\un$ by $\un_n:=\prodint{\un,x^n}, n\in\N$. For $c\in\mathbb{C}$ and $m\in\N$ we define the linear functionals $(x-c)^m\un$ and $(x-c)^{-m}\un$ by
\begin{equation}\label{deri}
\prodint{(x-c)^m\un,p(x)}=\prodint{\un,(x-c)^mp(x)}, \quad p\in\mathbb{P},$$ and $$ \prodint{(x-c)^{-m}\un,p(x)}=\prodint{\un,\dfrac{p(x)-\sum\limits_{k=0}^{m-1}\dfrac{D^{k}p(c)}{k!}(x-c)^k}{(x-c)^m}},\quad p\in\mathbb{P},
\end{equation}
where $D$ denotes the usual derivative.

\begin{defi}
Let $\un$ and $\vn$ be two linear functionals.
\begin{enumerate}
	\item The derivative of $\un$   is defined as
	$$\prodint{\un^\prime,p(x)}=-\prodint{\un,p^\prime(x)}.$$
	\item The sum and the  product of $\un$ and $\mathbf{v}$ are defined from their moments as follows \cite{Maro1991}
	$$(\un +\vn)_{n}=\un_n+\vn_n.$$
	$$ (\un \vn)_{n}= \prodint{\un\mathbf{v},x^n}=\sum_{k=0}^n\un_k\vn_{n-k},\quad n\geq 0.$$
	This product is commutative, associative and distributive with respect to the sum of linear functionals.
\end{enumerate}
\end{defi}

Let $c$ be a complex number and let $\deltan_c$ be the linear functional defined by $$\prodint{\deltan_c,x^n}=c^n,\quad  n\in\N. $$ It is not difficult to check that for any linear functional $\un,$ $\un\deltan_0=\un$. Moreover, if the first moment of $\un$ is nonzero, then there exists a unique linear functional $\un^{-1}$ such that $\un\un^{-1}=\deltan_0$. The moments of $\un^{-1}$ are defined recursively by
\begin{equation}\label{invermoments}
(\un^{-1})_{n}=-\frac{1}{\un_0}\sum_{k=0}^{n-1}\un_{n-k}(\un^{-1})_k,\quad n\geq 1,\quad  (\un^{-1})_{0}=\un_0^{-1}.
\end{equation}
The linear functional $\un$ is said to be quasi-definite when every leading principal submatrix of the Hankel matrix $H=(\un_{i+j})_{i , j=0}^{\infty}$ is nonsingular.  In such a situation, there exists a sequence of monic polynomials $(P_n)_{n\geq 0}$ such that $\deg{P_n}=n$ and  $\prodint{\un,P_n(x)P_m(x)}=K_{n}\delta_{n,m},$ where $\delta_{n,m}$ is the Kronecker  symbol and $K_{n}\ne 0$ (see \cite{Chihara}). The sequence $(P_n)_{n\geq 0}$ is said to be the sequence of monic orthogonal polynomials (SMOP) with respect to the linear functional $\un$. If every leading principal submatrix of the Hankel matrix is positive definite, then $\un$ is said to be a positive-definite linear functional. \\

If $\un$ is a quasi-definite linear functional and $(P_n)_{n\geq 0}$  is its corresponding  SMOP,  then there  exist two  sequences of complex numbers $( a_n)_{n\geq 1}$  and   $(b_n)_{n\geq 0}$, with $ a_n\ne 0$, such that	
\begin{equation}\label{ttrr}
\begin{aligned}
x\,P_{n}(x)&=P_{n+1}(x)+b_n\,P_{n}(x)+ a_{n}\,P_{n-1}(x),\quad n\ge 0,\\
P_{-1}(x)&=0, \ \ \ \ \ P_{0}(x)=1.
\end{aligned}
\end{equation}
Conversely, from Favard’s Theorem (see \cite{Chihara}) if $(P_n)_{n\geq 0}$ is a sequence of monic polynomials  generated by a three term recurrence relation  as in \eqref{ttrr} with $ a_n\ne 0,$ $n\geq 1$, then there exists a unique linear functional $\un$ such that $(P_n)_{n\geq0}$ is its SMOP.


Another way to write the recurrence relation \eqref{ttrr} is in a matrix form. Indeed, if ${\bf P}=(P_0,P_1, \cdots)^\top,$  where $A^\top$ denotes the transposed of a matrix $A,$ then   $x{\bf P}=J{\bf P},$ where $J$ is the  semi-infinite matrix
\begin{equation*}
J=\begin{pmatrix}
b_0 &1  &   &\\
a_1 &b_1& 1&\\
& a_2&b_2&\ddots\\
& & \ddots&\ddots
\end{pmatrix}.
\end{equation*}
The matrix $J$ is known in the literature as monic Jacobi matrix (see \cite{Chihara}).
\begin{defi}
Let $(P_n)_{n\geq 0}$ be the SMOP with respect to the linear functional  $\un$  satisfying the three term recurrence relation~\eqref{ttrr}.  For $k\in\N$ we define the sequence of associated polynomials of the $k$-th kind $(P^{(k)}_n)_{n\geq 0}$,  (also called the $k$-th associated polynomials, see \cite{Chihara})   as the sequence of monic polynomials satisfying the following recurrence relation
\begin{align}\label{asociados}
xP^{(k)}_n(x)&=P^{(k)}_{n+1}(x)+b_{n+k}P^{(k)}_n(x)+ a_{n+k}P^{(k)}_{n-1}(x),  \quad n\ge 0, \\ P^{(k)}_{-1}(x)&=0, \quad P^{(k)}_{0}(x)=1.\notag
\end{align}
\end{defi}
According to Favard's Theorem, there exists a quasi-definite linear functional $\un^{(k)}$, called the $k$-associated transformation of $\un$, such that $(P_n^{(k)})_{n\geq 0}$ is  its corresponding  SMOP.
The associated polynomials of the $k$-th kind $(P_n^{(k)})_{n\geq 0}$  can be expressed  as
\begin{equation}\label{asociadosrepr}
P_{n-1}^{(k)}(x)=\frac{1}{\un^{(k-1)}_0}\prodint{\un^{(k-1)}_y,\frac{P^{(k-1)}_n(x)-P^{(k-1)}_n(y)}{x-y}},\quad n\ge 1.
\end{equation}

Here $\un^{(k)}_y$ means that the linear functional $\un^{(k)}$ acts on the variable $y$.

On the other hand, there is a direct representation of such polynomials as (see \cite{Belme, Walter})
$$P_{n-k}^{(k)}(x)=\frac{1}{\prodint{\un, P_{k-1}^{2}}}\prodint{P_{k-1}(y)\un_y,\frac{P_n(x)-P_n(y)}{x-y}},\quad n\ge k. $$
Taking into account that $(P_n(x))_{n\geq 0}$ and $(P_{n-1}^{(1)}(x))_{n\geq 0}$ are two linearly independent solutions of  the difference equation  (see \cite{Walter})
$$xs_n=s_{n+1}+b_ns_n+ a_ns_{n-1},\quad n\geq 1,$$
every solution can be represented as a linear combination of $(P_n(x))_{n\ge 0}$ and $(P_{n-1}^{(1)}(x))_{n\ge 0}$. In particular (see \cite{Marcellan,Walter}),
\begin{equation}\label{linearcombinationk}
P_{n-k}^{(k)}(x)=A(x,k)P_n(x)+B(x,k)P_{n-1}^{(1)}(x), \quad n\ge k,
\end{equation}
where
$$A(x,k)=-\dfrac{P^{(1)}_{k-2}(x)}{\prod_{m=1}^{k-1} a_m}\quad \text{and} \quad B(x,k)=\dfrac{P_{k-1}(x)}{\prod_{m=1}^{k-1} a_m}.$$

Given a quasi-definite linear functional $\un$, we can define the formal series $$ \Su_{\un}(z):=\sum_{n=0}^\infty \frac{\un_n}{z^{n+1}},$$  that is said to be the Stieltjes   function  associated with  $\un$.  Using the coefficients of the three term recurrence relation \eqref{ttrr} we can represent the Stieltjes function in terms of a continued fraction (see~\cite{Chihara})
\begin{align}
\Su_{\un}(z)&=\cfrac{\un_0}
{(z-b_0) -\cfrac{ a_1}{(z-b_1)-\cfrac{ a_2}{(z-b_2)-\cfrac{ a_3}{(z-b_3)-\ddots}}}
}\label{stieltjes1}\\
&= \cfrac{\un_0}{(z-b_0)-\cfrac{ a_1}{\un^{(1)}_0} \Su_{\un^{(1)}}(z)},\label{stieltjes}	
\end{align}

where $\Su_{\un^{(1)}}(z)$ is the Stieltjes   function associated with $\un^{(1)}$. In order to have rational approximations of the Stieltjes function $\Su_{\un}(z)$ the standard way is to use the diagonal Padé approximants \cite{Pade,Stahl97}. Indeed,
$$	 \Su_{\un}(z)=\frac{\un_0P_{n-1}^{(1)}(z)}{P_n(z)}+\mathcal{O}\left(\frac{1}{z^{2n+1}}\right).$$

From here the associated polynomials of the first kind are also known in the literature as numerator polynomials (see \cite{Chihara}).

If $\un$ is quasi-definite, \eqref{invermoments} yields \begin{equation}\label{identidad}
\Su_{\un}(z) \Su_{\un^{-1}}(z)=1/z^2\end{equation} and this can be used together  with \eqref{stieltjes} to get (see \cite{Alf2004,Maro1988})

\begin{equation}\label{relationS}
\Su_{\un^{(1)}}(z)=-\frac{\un_0\un^{(1)}_0}{ a_1}z^2 \Su_{\un^{-1}}(z)+\dfrac{\un^{(1)}_0}{ a_1}(z-b_0).
\end{equation}
Hence
\begin{equation}\label{fu1}
\un^{(1)}= - \dfrac{\un_0^{(1)}\un_0}{ a_1}x^2\un^{-1}.
\end{equation}

In other words, $\un^{-1}$ is a quadratic Geronimus transformation of $\un^{(1)}$ (see~\cite{AMPR11}).

\begin{defi}
Let $(P_n)_{n\geq 0}$ be a SMOP with respect to $\un$  satisfying the recurrence relation \eqref{ttrr}. The sequence of monic polynomials $(P_n(x;\alpha))_{n\geq 0}$ is said to be  co-recursive  of parameter $\alpha$ associated with the linear functional $\un$,  if they also  satisfy \eqref{ttrr} but  with initial conditions $P_0(x;\alpha)=1$, and $P_{1}(x;\alpha)=P_1(x)-\alpha$. Notice that $$P_n(x;\alpha)=P_n(x)-\alpha P^{(1)}_{n-1}(x),\quad n\geq 0.$$
\end{defi}
Using the above relation we can check that the polynomials $P_{n}(x,\alpha)$ satisfy the recurrence relation
\begin{equation}\label{ttrr1}
\begin{aligned}
&x\,P_{n}(x;\alpha)=P_{n+1}(x;\alpha)+b_n\,P_{n}(x;\alpha)+ a_{n}\,P_{n-1}(x;\alpha),\quad n\ge 1,\\
&x\,P_{0}(x;\alpha)=P_{1}(x;\alpha)+(b_0+\alpha)\,P_{0}(x;\alpha).\\
\end{aligned}
\end{equation}

Let $\un^\alpha$ be the linear functional associated with the sequence of co-recursive polynomials $(P_n(x;\alpha))_{n\geq 0}$, then from    \eqref{stieltjes1} and \eqref{ttrr1}
$$\Su_{\un^\alpha}(z)=\cfrac{(\un^\alpha)_0}{(z-b_0-\alpha)-\cfrac{ a_1}{\un^{(1)}_0} \Su_{\un^{(1)}}(z)},$$ or, equivalently,

\begin{equation*}
\frac{1}{\un^\alpha_0}\left((z-b_0-\alpha)-\cfrac{ a_1}{\un^{(1)}_0} \Su_{\un^{(1)}}(z)\right) \Su_{\un^\alpha}(z)=1.\end{equation*}

Introducing \eqref{relationS} in the above formula
$$ \Su_{\un^\alpha}(z)\Biggl[\frac{-\alpha}{\un^\alpha_0z^2}+\frac{\un_0}{\un^\alpha_0} \Su_{\un^{-1}}(z) \Biggr]=\frac{1}{z^2}.$$
Thus, from \eqref{identidad} we get that $\dfrac{-\alpha}{\un^\alpha_0z^2}+\dfrac{\un_0}{\un^\alpha_0} \Su_{\un^{-1}}(z) $ is the Stieltjes function associated with $(\un^{\alpha})^{-1}.$ Comparing their moments we obtain
$(\un^{\alpha})^{-1}=\dfrac{\un_0}{\un^\alpha_0}\un^{-1}+\dfrac{\alpha}{\un^\alpha_0}\deltan^\prime_0$.  From here
\begin{equation}\label{funccorre}
\un^{\alpha}=\dfrac{\un_0^\alpha}{\un_0}\left(\un^{-1}+\dfrac{\alpha}{\un_0}\deltan^\prime_0\right)^{-1}.
\end{equation}
The above result is not new (see \cite{Maro1991}), however this alternative proof shows the way how you can find from the Stieltjes  an expression of the linear functional $(\un^{\alpha})^{-1}$  in terms of the linear functional $\un.$.\\

The analysis of perturbations of linear functionals is an interesting topic in the theory of orthogonal polynomials on the real line (scalar OPRL) (\cite{Bue1,Bue2,DereM,Dere,Zhed97} and references therein). In particular, when dealing with a positive definite case, i.e., the linear functional has an integral representation in terms of a probability measure supported in an infinite subset of the real line, such perturbations provide a useful information in the study of Gaussian quadrature rules for the perturbed linear functional, taking into account that the perturbation yields new nodes and Christoffel numbers \cite{Chihara,Gaut3}.

Among the perturbations of linear functionals,  spectral linear perturbations have attracted the interest of  researchers (see \cite{Zhed97}). Such perturbations are generated by two particular families, the so called Christoffel and Geronimus transformations.\\
Christoffel perturbations, that  appear when considering orthogonality with respect to a new linear functional $\widetilde{\mathbf{u}}= p(x) \mathbf{u}$, where $p(x)$ is a polynomial, were studied in 1858 by  E. B.  Christoffel  in the framework of Gaussian quadrature rules \cite{chr1877}. He  found  explicit connection formulas between  the corresponding sequences of orthogonal polynomials with respect to both measures,  the Lebesgue measure  $d\mu$ supported on the interval $(-1,1)$ and  $d\widetilde{\mu}(x)= p(x) d\mu(x)$, with $p(x)=(x-q_1)\cdots(x-q_N)$, a signed polynomial in the support of $d\mu$,  as well as the distribution of their zeros as nodes in such quadrature rules. Nowadays, these are  known in the literature as  Christoffel formulas (see \cite{Chihara,Sz78}). More recently in \cite{leary} the authors studied the sensitivity of Gauss–Christoffel quadrature with respect to small perturbations of the probability measure.  With regard to the corresponding sequence of monic orthogonal polynomials (SMOP in short),   explicit relations between the polynomials and the coefficients of the three term recurrence relations that they satisfy  have been extensively studied, see \cite{Gau82}, as well as the relation between  the Jacobi matrices in the framework of the so-called discrete Darboux transformations. They are  based on the $LU$ factorization of such matrices (see \cite{Bue1} and \cite{Yoon02}, among others).

%
%
It is worthwhile to note that the zeros of orthogonal polynomials with respect to the canonical Christoffel transformation of a positive definite linear functional are the nodes in the Gauss-Radau quadrature formula. In the case of a perturbation of the measure by a positive quadratic polynomial in the support of the measure, the zeros of the corresponding orthogonal polynomials are the nodes of the Gauss-Lobatto quadrature rule (see \cite{Gaut02}).

Geronimus transformations appear when you deal with perturbed functionals $ \widehat \un$ defined by $p(x) \widehat \un=\mathbf{u}$, where $p(x)$ is a polynomial. Such a kind of transformations were used by  J. L. Geronimus (see \cite{Ger40}), in order to provide an alternative proof of a result given by W. Hahn \cite{Hahn35} concerning the characterization of classical orthogonal polynomials (Hermite, Laguerre, Jacobi and Bessel) as the unique families of orthogonal polynomials whose first derivatives are also orthogonal polynomials.
%
Examples of such transformations have been given by P. Maroni \cite{Maro1990} for a perturbation of the type $p(x)=x-c, $ in \cite{AMPR11,BM98,Buh} for a quadratic case, and in \cite{MN03} for the cubic case.

As it was mentioned above the Christoffel and Geronimus transformation are known in the literature as discrete Darboux transformations due to the existing relation between their Jacobi matrices associated  and  LU and UL factorizations~\cite{Bue1,Yoon02}.
\begin{equation*}
\begin{array}{|c|c|}
\hline
\begin{array}{c|c}
\un &\widetilde\un=(x-c)\un\\
J-cI=LU&\widetilde J-cI=UL
\end{array}& \substack{\text{Darboux transformation}\\ \text{without parameter}}\\
\hline
\begin{array}{c|c}
\un&\un=(x-c)\widehat\un\\
J-cI=UL&\widehat J-cI=LU\end{array}& \substack{\text{Darboux transformation}\\ \text{}}\\
\hline
\end{array}
\end{equation*}

They correspond to the factorization of a discrete operator of second order in terms of two discrete operators of first order. We have focused the attention on these transformations because they constitute a generating system of the set of linear spectral transformations, i.e., every linear spectral transformation can be presented as a finite composition of Christoffel and Geronimus transformation~(see \cite{Zhed97}). \\

With the above background we are ready to point out the goals of our contribution.  Let $\un$ be a quasi-definite functional and let $\tilde\un=(x-c)\un$ and $(x-c)\hat \un=\un$ be the canonical Christoffel and Geronimus transformation of $\un,$ respectively. As illustrated the Figure \ref{fig:M1}, we are interested in deducing relations between $\un^{(1)}$ and $\widetilde\un^{(1)}$ (resp. $\widehat\un^{(1)}$)  by using   the LU and UL factorization of the monic Jacobi matrix associated with $\un$, as well as explicit algebraic relations between their corresponding  SMOP. With this in mind the structure of the manuscript is as follows.

{ \begin{figure}[h]
	\huge\centering
	\begin{tikzcd}
	\un \arrow{r}[near start]{\substack{\text{\small\quad\  associated}\\ \text{ \small \quad functional}}} \arrow[d, "\substack{\text{\small Darboux}\\\text{ \small transformation}}" left]
	& \un^{(1)} \arrow[d,blue,dashrightarrow, "{\bf ?}" black] \\
	\widetilde\un \arrow[r]
	& \widetilde{\un}^{(1)}
	\end{tikzcd}
	\caption{Structure of the problem.} \label{fig:M1}
\end{figure}
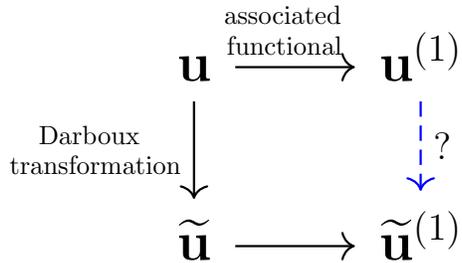}

In Section \ref{section2} we study a quadratic Geronimus transformation of a linear functional assuming that the zero of the polynomial has multiplicity equal to 2. The relation between the coefficients of the three term recurrence relation that the corresponding SMOP satisfies is given as well as the connection between the respective Jacobi matrices. As an application, the coefficients of the three term recurrence relation when you deal with the inverse of a linear functional are deduced taking into account such a linear functional is a quadratic Geronimus transformation of the associated linear functional of the first kind. In Section \ref{section3}
we deal with a Darboux transformation (Christoffel or Geronimus) followed by an associated transformation of the first kind of the linear functional. In this way (see Figure \ref{fig:M1}) we focus the attention on the relation between the resulting linear functionals, their corresponding SMOP and their Jacobi matrices.

Specifically, we analyze the behavior of associated SMOP of the first kind under canonical Christoffel (and Geronimus) transformations.  Finally, in Section~\ref{section4} we show some examples related to the application of these linear transformations for classical linear functionals and their corresponding associated linear functionals of the first kind.

\section{Orthogonal polynomials with respect to the inverse of a linear functional} \label{section2}

\begin{pro}
Let $\un$ be a linear functional and let $\widehat\un$ be the linear functional defined by the following Geronimus transformation of $\un$:
\begin{equation}\label{cuadratic}
\un=(x-c)^2\widehat\un,\quad c\in\mathbb{C}.
\end{equation}

Then
$$  \Su_{\widehat\un} (z)=\frac{ \Su_{\un}(z)+A(z-c)+B}{(z-c)^2},$$
where $A=\widehat{\un}_0$ and $B=\widehat{\un}_1-c\widehat{\un}_0$.
\end{pro}

\begin{proof}
Although the proof of this result can be obtained easily from the results given in \cite{Bue1,DereM,Dere} here we give a complete proof.  	

Computing the moments of $\un$	
$$\prodint{\un, x^n} =\prodint{(x-c)^2 \widehat\un, x^{n}}= \widehat\un_{n+2} - 2c \widehat\un_{n+1} + c^{2} \widehat\un_{n}, \quad n\geq 0.$$
Thus, the corresponding Stieltjes functions are related by
\begin{equation*}
\Su_{\un} (z)= z^{2}\left( \Su_{\widehat{\un}} (z) - \frac{\widehat\un_{0}} {z} - \frac{\widehat\un_{1}}{z^2}\right)- 2 c z \left( \Su_{\widehat\un} (z) - \frac{\widehat\un_{0}}{z}\right) + c^{2}  \Su_{\widehat\un}(z).
\end{equation*}
In other words,
$$  \Su_{\widehat\un} (z)=\frac{ \Su_{\un}(z)+A(z-c)+B}{(z-c)^2},$$
where $A=\widehat{\un}_0$ and $B=\widehat{\un}_1-c\widehat{\un}_0$.
\end{proof}
Computing the moment of $\widehat\un$ we get
\begin{align*}
\prodint{\widehat\un,x^n}&=\prodint{(x-c)^2\widehat\un,\dfrac{x^n-c^n-nc^{n-1}(x-c)}{(x-c)^2}}+c^n\prodint{\widehat\un,1}+nc^{n-1}\prodint{\widehat\un,(x-c)}\\
&=\prodint{(x-c)^{-2}\un,x^n}+\widehat\un_0\prodint{\deltan_c,x^n}-[\widehat\un_1-c\widehat\un_0]\prodint{\deltan^\prime_c,x^n}.
\end{align*}
In terms of linear functionals, the above relation reads (see \cite{Zhed97})
\begin{equation}\label{repre}
\widehat\un = (x-c)^{-2} \un + \widehat\un_{0} \deltan_{c} - \left[\widehat\un_{1} - c \widehat\un_{0}\right] \deltan'_{c}.
\end{equation}

Notice that in the definition of   $\widehat\un$ we have two degrees of freedom, corresponding to the choices of  $\widehat\un_{0}$, and~$ \widehat\un_{1}.$

\begin{pro}Let $(P_n)_{n\geq 0}$ be the SMOP associated with the linear functional $\un$, then
\begin{equation}\label{deriv}
(P_{n-1} ^{(1)})' (c)=\frac{1}{\un_{0}} \prodint{(x-c)^{-2} \un, P_{n}(x)}.
\end{equation}

\end{pro}
\begin{proof}
Since $P_{n-1}^{(1)}(x)= \dfrac{1}{\un_{0}}\prodint{\un_{y}, \frac{ P_{n}(x)- P_{n}(y)}{x-y}}$, taking derivatives with respect to the variable $x$ in the above expression and evaluating  it at $ x=c$ the result follows immediately.
\end{proof}

\begin{pro}\label{pro4} If $(P_n)_{n\geq 0}$ is the SMOP associated with the linear functional $\un$, then the linear functional $\widehat\un$ defined by the expression  $\un=(x-c)^2\widehat\un$ is quasi-definite if  and only if $\widehat\un_0\ne0$ and $d_{n}^{*}\neq0$, where
$$d_n^*:=\det\begin{pmatrix}
S_{n-2}(c)&S_{n-1}(c)\\
S_{n-2}^\prime(c)+
\widehat \un_0P_{n-2}(c)&S_{n-1}^\prime(c)+\widehat \un_0P_{n-1}(c)
\end{pmatrix},\quad n\geq 2.$$
Here $S_n(x)=[\widehat\un_{1} - c \widehat\un_{0}]P_n(x)+\un_0P^{(1)}_{n-1}(x)$.
Moreover, if $(Q_n)_{n\geq 0}$ is the   SMOP associated with $\widehat\un,$ then
\begin{align}
&Q_{n}(x)=\frac{1}{d_n^*}\times \label{repreQ}\\
&\det\begin{pmatrix}
P_n(x)&P_{n-1}(x)&P_{n-2}(x)\\
S_{n}^\prime(c)+\widehat \un_0P_{n}(c)&S_{n-1}^\prime(c)+\widehat \un_0P_{n-1}(c)&S_{n-2}^\prime(c)+\widehat \un_0P_{n-2}(c)\\
S_{n}(c)&S_{n-1}(c)&S_{n-2}(c) 	
\end{pmatrix},\  n\geq 2\notag,
\end{align}
$$Q_0(x)=1,\quad Q_{1} (x)= x - \frac{\widehat\un_{1}}{\widehat\un_{0}}. $$
\end{pro}
\begin{proof}
Assuming $\widehat\un$ is quasi-definite, let us consider the Fourier expansion
\begin{equation}\label{re}
Q_{n}(x)= P_{n}(x)+\sum_{m=0}^{n-1}\alpha_{n,m}P_{m}(x). \end{equation}
Here
$$ \alpha_{n,m}=\frac{\prodint{\un ,  Q_{n}(x) P_{m}(x)}}{ \prodint{\un ,  P^{2}_{m}(x)}}= 0, \quad m\le n-3,$$
$$\alpha_{n,n-2}= \frac{\prodint{\un,  Q_{n}(x) P_{n-2}(x)}}{\prodint{\un, P^{2}_{n-2}(x)}}= \frac{\prodint{\widehat\un,  Q^{2}_{n}(x)}}{\prodint{\un, P^{2}_{n-2}(x)}} \neq 0.$$
Now, from \eqref{re}, we have the system, $n\geq 2$
\begin{equation}\label{sys1}
\Scale[0.98]{-\begin{pmatrix}
	\prodint{\widehat\un,P_{n}}\\
	\prodint{\widehat\un,(x-c)P_{n}}
	\end{pmatrix}=
	\begin{pmatrix}
	\prodint{\widehat\un,P_{n-2}}&\prodint{\widehat\un,P_{n-1}}\\
	\prodint{\widehat\un,(x-c)P_{n-2}}&\prodint{\widehat\un  ,(x-c)P_{n-1}}
	\end{pmatrix}\begin{pmatrix}
	\alpha_{n,n-2}\\
	\alpha_{n,n-1}
	\end{pmatrix}.}
\end{equation}
Taking into account \eqref{repre} and  \eqref{deriv} we get
\begin{align*}
&\prodint{\widehat\un,P_{n}}=\un_0 (P_{n-1} ^{(1)})^\prime (c)+\widehat\un_0P_n(c)+[\widehat\un_{1} - c \widehat\un_{0}]P_n^\prime(c)=S_{n}^\prime(c)+\widehat \un_0P_{n}(c),\\
&\prodint{\widehat\un,(x-c)P_{n}}=\un_0 P_{n-1} ^{(1)}(c) +[\widehat\un_{1} - c \widehat\un_{0}]P_n(c)=S_n(c).
\end{align*}
Since $Q_n(x)$ is a monic polynomial of degree $n$, we know that \eqref{sys1} has at least one solution.
Otherwise, if  we suppose that it has two different solutions, then there are two monic polynomials of degree $n$ satisfying the orthogonality condition. But this contradicts the uniqueness of the sequence $(Q_n) _ {n\geq0}.$    Thus $d_n^*\ne0$.

Conversely, let assume that $d_n^*\ne 0$ and define the polynomials $Q_{0}(x)=1,$ $Q_{1}(x)= x - \frac{\widehat\un_{1}}{\widehat\un_{0}}$ and
$Q_{n}(x)$ as in \eqref{repreQ} for $n\geq 2$.
Then it is not difficult to check that $(Q_n)_{n\geq 0}$ is the SMOP  with respect to the functional $\widehat\un.$ Notice that the representation of   $Q_n(x)$ depends on the first two  moments $\widehat\un_{0}$ and $ \widehat\un_{1}.$
\end{proof}
\begin{rem}
Observe that for $n\geq 2,$
\begin{equation}\label{eq10}
\begin{aligned}
\alpha_{n,n-1}&=-\dfrac{1}{d_n^*}\det\begin{pmatrix}
S_{n}^\prime(c)+
\widehat \un_0P_{n}(c)&S_{n-2}^\prime(c)+\widehat \un_0P_{n-2}(c)\\
S_{n}(c)&S_{n-2}(c)\\
\end{pmatrix},\\
\alpha_{n,n-2}&=\dfrac{d^*_{n+1}}{d^*_{n}}.	
\end{aligned}
\end{equation}

\end{rem}

When $\widehat\un$ is quasi-definite, the coefficients of the three term recurrence relation associated with $(Q_n)_{n\geq 0}$ can be written as follows (see \cite{AMPR11} for a alternative proof).
\begin{pro}\label{ttrrq}
Let $(P_n)_{n\geq 0}$ be the SMOP with respect to $\un$  satisfying the three term recurrence relation \eqref{ttrr}. If $\widehat\un$ is quasi definite, then
\begin{equation*}\label{1.2}
\begin{aligned}
&x\,Q_{n}(x)=Q_{n+1}(x)+\widehat b_n\,Q_{n}(x)+\widehat a_{n}\,Q_{n-1}(x),\quad n\ge 0,\\
&Q_{-1}(x)=0, \ \ \ \ \ Q_{0}(x)=1,
\end{aligned}
\end{equation*}
with
$$\begin{array}{l}{\widehat{b}_{n}=b_{n}+\alpha_{n,n-1}-\alpha_{n+1,n}, \  n \geq 0,\quad \widehat{ a}_{n}=\dfrac{\alpha_{n,n-2}}{\alpha_{n-1,n-3}} a_{n-2}, \  n \geq 3,}\\[10pt] { \widehat a_2=\dfrac{\un_0\widehat\un_0\alpha_{2,0}}{\un_0\widehat\un_0-(\widehat\un_{1}-c\widehat\un_0)^2},\quad \widehat a_1=\dfrac{\un_0\widehat\un_0-(\widehat\un_{1}-c\widehat\un_0)^2}{\widehat\un_0^2}}.\end{array}$$
\end{pro}
\begin{proof}
Since
$$\widehat a_n=\dfrac{\prodint{\widehat\un,(x-c)^nQ_n(x)}}{\prodint{\widehat \un,(x-c)^{n-1}Q_{n-1}(x)}},\quad n\geq 1,\quad \widehat b_n=\dfrac{\prodint{\widehat \un,xQ^2_n(x)}}{\prodint{\widehat \un,Q^2_n(x)}} ,\quad n\geq 0,$$	
the proof is a direct consequence of (\ref{repre}) \eqref{re} and \eqref{ttrr}.
\end{proof}

Now, taking into account that
$$\prodint{\widehat\un,(x-c)^2P_{n}(x)Q_k(x)}=\prodint{\un,P_{n}(x)Q_k(x)}=0,\quad k=0,\ldots,n-1,$$ there exist complex numbers $\beta_{n,n+1}$ and $\beta_{n,n}\ne 0$ such that
\begin{equation}\label{conex2}
(x-c)^2P_n(x)=Q_{n+2}(x)+\beta_{n,n+1}Q_{n+1}(x)+\beta_{n,n}Q_{n}(x), \quad n\geq0.
\end{equation}
Taking derivatives in (\ref{conex2}) and evaluating in  $x=c$ we get
\begin{equation}\begin{array}{cl}
\begin{cases}\label{system}
-Q_{n+2}(c)=Q_{n+1}(c)\beta_{n,n+1}+Q_n(c)\beta_{n,n},\\
-Q_{n+2}'(c)=Q_{n+1}'(c)\beta_{n,n+1}+Q_n'(c)\beta_{n,n},
\end{cases}&\quad n\geq0.\end{array}\end{equation}

Since $(Q_n)_{n\geq0}$ is a SMOP, the representation \eqref{conex2} is unique. This implies taht the system \eqref{system} has a unique solution and therefore   $$W(Q_{n+1},Q_n)(c)\neq0,$$ where
\begin{equation*}\label{wrosk}
W(p,q)(x)=p(x)q^\prime(x)-p^\prime(x)q(x), \quad p,q\in \mathbb{P}.
\end{equation*}
As a consequence,  we can give an explicit representation of  $\beta_{n,n}\ne0$ and $\beta_{n,n+1}$ as follows
\begin{equation}\label{betas}
\beta_{n,n}=\dfrac{W(Q_{n+1},Q_{n+2})(c)}{W(Q_{n},Q_{n+1})(c)}\quad \text{and}\quad \beta_{n,n+1}=-\dfrac{W(Q_{n},Q_{n+2})(c)}{W(Q_{n},Q_{n+1})(c)}.
\end{equation}

Notice that \eqref{re} and \eqref{conex2} can be written in a matrix form as
$$\mathbf{Q}=L\mathbf{P},\quad (x-c)^2\mathbf{P}=U\mathbf{Q},$$
where $\mathbf{P}=(P_0(x),P_1(x),\cdots)^\top,$  $\mathbf{Q}=(Q_0(x),Q_1(x),\cdots)^\top$ and
\begin{equation}\label{LUinversa}
\Scale[0.94]{L=\begin{pmatrix}
1&&&\\
\alpha_{1,0}&1&&\\
\alpha_{2,0}&\alpha_{2,1}&1&&\\
\vdots          & \ddots     &\ddots &\ddots
\end{pmatrix},\
U=\begin{pmatrix}
\beta_{0,0}& \beta_{0,1}&1&&&\\
& \beta_{1,1}& \beta_{1,2}&1&&\\
&& \beta_{2,2}& \beta_{2,3}&1&&\\
&&&\ddots&\ddots&\ddots
\end{pmatrix},}
\end{equation}
with $\alpha_{n,n-1}$, $\alpha_{n,n-2}$, $\beta_{n,n}$ and $\beta_{n,n+1}$  given by (\ref{eq10}) and (\ref{betas}), respectively. The above result is summarized in the following
\begin{pro}\label{propLUinversa}
Assume that $\un$ and $\widehat\un$ are quasi-definite linear functionals and let $J $ and $\widehat J $ be the corresponding monic Jacobi matrices. Then,
$$\left(J -cI\right)^2=UL \quad \text{and}\quad \left(\widehat J -cI\right)^2=LU,$$
where $L$ and $U$ are defined in \eqref{LUinversa}.
\end{pro}

Next, let us assume $\un$ and $\un^{-1}$ are  quasi-definite linear functionals and  let $(P^-_n)_{n\geq 0}$ be the SMOP with respect to  $\un^{-1}$.
From \eqref{fu1} and \eqref{repre} we get an explicit expression for $\un^{-1}$.

\begin{equation}\label{uminus}
\un^{-1}=-\dfrac{ a_1}{\un^{(1)}_0\un_0}x^{-2}\un^{(1)}+\dfrac{1}{\un_0}\deltan_0+\dfrac{\un_1}{\un_0^2}\deltan_0^\prime.
\end{equation}

On the other hand, by  \eqref{fu1}
$$\prodint{\un^{-1},x^{m}P_n^-(x)}=\prodint{\un^{(1)},x^{m-2}P_n^-(x) }=0,\quad m=2,3,\ldots, n-1.$$
Thus, there exist complex numbers $\alpha_{n,n-1}$ and $\alpha_{n,n-2}$, $n\ge2,$   such that
\begin{equation*}\label{conex}
P_n^{-}(x)=P_{n}^{(1)}(x)+\alpha_{n,n-1}P_{n-1}^{(1)}(x)+\alpha_{n,n-2}P_{n-2}^{(1)}(x), \quad n\ge 2,
\end{equation*}
$$P_{1}^{-} (x)= x + b_{0}.$$
With this in mind and taking into account Proposition \ref{pro4} and Proposition \ref{ttrrq}    we can state the following two useful Corollaries. See also~\cite{Maro1991}.
\begin{coro}\label{pnminusexpre}
If $\un$ is quasi-definite, then $\un^{-1}$ is  quasi-definite if and only if $b_1+b_0\ne 0$ and  $d_n^*\ne0$, $n\geq 2$, where
$$d_n^*=\dfrac{1}{\un_0^2}W(P_n,P_{n-1})(0).$$
Moreover, if for all $n\geq 2,$ $d_n^*\ne 0$, then
\begin{equation}\label{inv}
P^-_{n}(x)=\frac{1}{\un_0^2d_n^*}\det\begin{pmatrix}
P_n^{(1)}(x)&P^{(1)}_{n-1}(x)&P^{(1)}_{n-2}(x)\\
P_{n+1}(0)&P_{n}(0)&P_{n-1}(0)\\
P^\prime_{n+1}(0)&P^\prime_{n}(0)&P^\prime_{n-1}(0)
\end{pmatrix},\quad n\geq 2,\end{equation}
with $P^{-}_{0}(x)=1,$ and $P^{-}_{1}(x)=x+b_0.$
\end{coro}
\begin{proof}
Assume that  $\un^{-1}$ is quasi-definite  and let us define $\vn=-\dfrac{ a_{1}}{\un_0\un_0^{(1)}}\un^{(1)}$ and $\widehat\vn=\un^{-1},$ respectively. Taking into account that $\vn=x^2\widehat \vn$  and   $\un^{(1)}$ is also quasi-definite, the proof is a direct consequence of   Proposition \ref{pro4} with $$S_n(x)=-\dfrac{\un_1}{\un^2_0}P_n^{(1)}(x)-\frac{ a_1}{\un_0}P_{n-1}^{(2)}(x).$$ Since for $n\geq 1$
\begin{align*}
S_{n}^\prime(0)+\widehat \vn_0P^{(1)}_{n}(0)=-\dfrac{ a_1}{\un_0}\left(P^{(2)}_{n-1}\right)^\prime(0)+\dfrac{1}{\un_0}P_n^{(1)}(0)-\dfrac{\un_1}{\un^2_0}\left(P_n^{(1)}\right)^\prime(0),
\end{align*}
taking into account that from \eqref{linearcombinationk}
$$\left(P^{(2)}_{n-1}\right)^\prime(0)=\frac{1}{ a_1}\left(-P^\prime_{n+1}(0)+P_n^{(1)}(0)-\dfrac{\un_1}{\un_0}\left(P_n^{(1)}\right)^\prime(0)\right),$$
we get
$$S_{n}^\prime(0)+\widehat \vn_0P^{(1)}_{n}(0)=\frac{1}{\un_0}P_{n+1}^\prime(0),\quad n\geq 1.$$
In the same way, keeping  in mind that
\begin{align*}
S_n(0)&=-\frac{ a_1}{\un_0}P_{n-1}^{(2)}(0) -\dfrac{\un_1}{\un^2_0}P_n^{(1)}(0).
\end{align*}
and
$$P_{n-1}^{(2)}(0)=\dfrac{1}{ a_1}\left(-P_{n+1}(0)-\dfrac{\un_1}{\un_0}P_n^{(1)}(0)\right),$$
then
$$S_n(0)=\dfrac{1}{\un_0}P_{n+1}(0),\quad n\geq 1.$$
\end{proof}
\begin{coro}\label{ttrrpminus}
The SMOP $(P^-_n)_{n\geq0}$ satisfies the following recurrence relation
$$xP^-_n(x)=P^-_{n+1}(x)+b_n^-P^-_n(x)+ a_n^-P^-_{n-1}(x),$$
where
\begin{align}\label{coe}
b^-_{n}&=b_{n+1}-\dfrac{W(P_{n+1},P_{n-1})(0)}{W(P_{n},P_{n-1})(0)}+\dfrac{W(P_{n+2},P_{n})(0)}{W(P_{n+1},P_{n})(0)},\quad n\geq 0,\\[10pt]
\notag a^-_{n}&=\dfrac{W(P_{n+1},P_{n})(0)W(P_{n-1},P_{n-2})(0)}{W(P_{n},P_{n-1})^2(0)} a_{n-1},\quad n\geq 2,\\
\notag a^-_1&=-(b_0^2+ a_1).
\end{align}

\end{coro}

\begin{rem}
Notice that from the expression for $ a^-_1$, the linear functionals $\un^{-1}$ and $\un$ cannot be  positive-definite simultaneously.
\end{rem}
In the same spirit of the proof of Proposition \ref{propLUinversa},
we have the next result.
\begin{pro}\label{relationlu}
Let $\un^{(1)}$ and $\un^{-1}$ be quasi-definite linear functionals and  let $J ^{(1)}$ and $J ^-$ be the  corresponding monic Jacobi matrices. Then,
$$\left(J^{(1)} \right)^2=UL \quad \text{and}\quad \left(J^{-} \right)^2=LU,$$
where
\begin{equation*}
\Scale[0.98]{	L=\begin{pmatrix}
	1&&&\\
	\alpha_{1,0}&1&&\\
	\alpha_{2,0}&\alpha_{2,1}&1&\\
	\vdots          & \ddots     &\ddots &\ddots
	\end{pmatrix},\quad
	U=\begin{pmatrix}
	\beta_{0,0}& \beta_{0,1}&1&&&\\
	& \beta_{1,1}& \beta_{1,2}&1&&\\
	&& \beta_{2,2}& \beta_{2,3}&1&&\\
	&&&\ddots&\ddots&\ddots
	\end{pmatrix},
}
\end{equation*}
and
\begin{align*}
\beta_{n,n}&=\dfrac{W(P^-_{n+1},P^-_{n+2})(0)}{W(P^-_{n},P^-_{n+1})(0)},& \beta_{n,n+1}&=-\dfrac{W(P^-_{n},P^-_{n+2})(0)}{W(P^-_{n},P^-_{n+1})(0)},\\[5pt]
\alpha_{n,n-1}&=-\dfrac{W(P_{n+1},P_{n-1})(0)}{W(P_{n},P_{n-1})(0)}, & \alpha_{n,n-2}&=\dfrac{W(P_{n+1},P_{n})(0)}{W(P_{n},P_{n-1})(0)}.
\end{align*}

\end{pro}
Another way to give a relation between the monic Jacobi matrices associated with two linear functionals related by \eqref{cuadratic} (when they are positive-definite) is through the QR factorization \cite{Buh}. Since in our case we known that $\un$ and $\un^{-1}$ cannot be  positive-definite linear functionals simultaneously this method cannot be used. However, the relation between those monic Jacobi matrices can be obtained by using the so-called hyperbolic factorization QR (see \cite{Bue2}). Let us consider the diagonal matrices $D_{p^-}$ and $D_{p^{(1)}}$ with $(D_{p^{-}})_{i,i}=\prodint{\un^{-1},(P^{-}_i)^2}$ and $(D_{p^{(1)}})_{i,i}=\prodint{\un^{(1)},\bigl(P^{(1)}_i\bigr)^2}$. These matrices do not have necessary positive main diagonal entries. However, introducing the diagonal matrices $|D_{p^{(1)}}|^{1/2}$ and $|D_{p^{-}}|^{1/2}$ such that $(|D_{p^{-}}|^{1/2})_{i,i}=|(D_{p^{-}})_{i,i}|^{1/2}$ and $(|D_{p^{(1)}}|^{1/2})_{i,i}=|(D_{p^{(1)}})_{i,i}|^{1/2}$ it turns out
$$D_{p^-}=|D_{p^{-}}|^{1/2}\Omega_{p^-}|D_{p^{-}}|^{1/2},\quad D_{p^{(1)}}=|D_{p^{(1)}}|^{1/2}\Omega_{p^{(1)}}|D_{p^{(1)}}|^{1/2},$$
where $\Omega_{p^-}$ and $\Omega_{p^{(1)}}$ are diagonal matrices with $(\Omega_{p^-})_{i,i}=\textnormal{sgn}\bigl(\prodint{\un^{-1},(P_i^{-})^2}\bigr)$ and $(\Omega_{p^{(1)}})_{i,i}=\textnormal{sgn}\bigl(\prodint{\un^{(1)},(P_i^{(1)})^2}\bigr)$. Finally, setting $\widetilde{D}_{p^-}=-\frac{\un_0^{(1)}\un_0}{ a_{1}}D_{p^-}$ it can be stated the following result.

\begin{pro}[\cite{Bue2}]
If all the leading principal submatrices of $J^- D_{p^-}(J^- )^\top $ are non-singular, then there exist  an upper triangular matrix $R$, with positive diagonal entries, and a matrix $Q$ with $Q^\top \Omega_{p^-}Q=\Omega_{p^{(1)}}$ such that $|\widetilde{D}_{p^-}|^{1/2}(J^- )^\top =QR$. Moreover, if $G=|D_{p^{(1)}}|^{1/2}Q^\top |\widetilde{D}_{p^-}|^{-1/2}$, then
$$J^- =LG,\quad J^{(1)} =GL,$$
where $L$ is given in Proposition \ref{relationlu}.
\end{pro}

\begin{proof}
Although the proof of this result can be obtained in the general framework given in \cite{Bue2} here we give a complete proof.

Since all the leading principal submatrices of $J^- D_{p^-}(J^- )^\top $ are non-singular there exist $R$, an upper triangular matrix with positive diagonal entries, and a matrix $Q$, with $Q^\top \Omega_{p^-}Q=\Omega_{p^{(1)}}$, such that (see \cite[Prop. 4.10]{Bue2}).	
\begin{equation}\label{QR}
|\widetilde{D}_{p^-}|^{1/2}(J^- )^\top =QR.
\end{equation}
On one hand, from \eqref{inv}, $\mathbf{P^-}=L\mathbf{P^{(1)}}$, where $L$ is given in Proposition \ref{relationlu}, $\mathbf{P^-}=(P_0^-(x),P_1^-(x),\ldots)^\top $ and $\mathbf{P^{(1)}}=(P_0^{(1)}(x),P_1^{(1)}(x),\ldots)^\top $. Thus we get
\begin{equation}\label{upp1}
\prodint{\un^{(1)},\mathbf{P^-}(\mathbf{P^-})^\top }=LD_{p^{(1)}}L^\top =L|{D}_{p^{(1)}}|^{1/2}\Omega_{p^{(1)}}|{D}_{p^{(1)}}|^{1/2}L^\top .
\end{equation}
On the other hand, using the definition of $J^{-} $ and \eqref{fu1} we also have
\begin{equation}\label{upp2}
\prodint{\un^{(1)},\mathbf{P^-}(\mathbf{P^-})^\top }=J^{-} \widetilde{D}_{p^-}(J^{-} )^\top =(QR)^\top \Omega_{p^-}Q R=R^\top \Omega_{p^{(1)}}R.
\end{equation}
From \eqref{upp1}, \eqref{upp2} and the uniqueness of the triangular factorization of the symmetric matrix $\prodint{\un^{(1)},\mathbf{P^-}(\mathbf{P^-})^\top }$ it turns out that
\begin{equation}\label{RL}
R=|D_{p^{(1)}}|^{1/2}L^\top .
\end{equation}
Thus, according to \eqref{QR} and \eqref{RL}
$$J^- =R^\top Q^\top |\widetilde{D}_{p^-}|^{-1/2}=LG,$$
where $G=|D_{p^{(1)}}|^{1/2}Q^\top |\widetilde{D}_{p^-}|^{-1/2}$. Finally, notice that
$\mathbf{P^{(1)}}=L^{-1}\mathbf{P^-}$, with the definition of $D_{p^{(1)}}$, yields  $J^{(1)} D_{p^{(1)}}=L^{-1}J^{-} LD_{p^{(1)}}$. As a consequence,
$$J^{(1)} =L^{-1}J^- L=GL.$$
\end{proof}
\section{ Associated orthogonal polynomials of first kind and Darboux transformations.} \label{section3}

\subsection{Christoffel transformation} Let $\un$ be a quasi-definite linear  functional and let $(P_n)_{n\geq 0}$ be its corresponding SMOP. If $c$ is a fixed complex number, the linear functional $\widetilde \un=(x-c)\un$ is said to be the canonical Christoffel transformation of the linear  functional $\un$. Suppose that $\widetilde\un$ is also quasi-definite (it is equivalent to $P_n(c)\ne 0$ for all $n\in\N$)
and let $(\widetilde P_n)_{n\geq 0}$ be its SMOP.  It is well known that $(P_n)_{n\geq 0}$ and $(\widetilde P_n)_{n\geq 0}$ are related by \cite{Chihara}
\begin{equation}\label{repChris}
(x-c)\widetilde{P}_n(x)=P_{n+1}(x)-\dfrac{P_{n+1}(c)}{P_n(c)}P_n(x), \quad n\geq 0.
\end{equation}
We have the following relation between their monic Jacobi matrices.
\begin{teo}[\cite{Bue1,Yoon02}]\label{teo1}
Let $J$ and $\widetilde J$ be the monic Jacobi matrices associated with $\un$ and  $\widetilde{\un} = (x-c)\,\un$, respectively. If $P_n(c)\ne 0$, for all $n\in\N$, then $J-cI$ has LU factorization, i.e.,
\begin{equation} \label{LU}
J-cI:=LU:=\begin{pmatrix}
1           &            &&\\
\ell_{1}&1           &&\\
&\ell_{2}&1&\\
&            &\ddots&\ddots
\end{pmatrix}\begin{pmatrix}
\beta_{0}&1            &             &&\\
&\beta_{1}&1            &&\\
&             &\beta_{2}&\ddots &\\
&             &             &\ddots&\ddots
\end{pmatrix},\end{equation}
where $L$ is a lower bidiagonal matrix with 1's as diagonal entries, $$\ell_n=-\frac{P_{n-1}(c)}{P_n(c)}\frac{\prodint{\un,P_n^2}}{\prodint{\un,P_{n-1}^2}},$$  and $U$ is an upper bidiagonal matrix with $\beta_n=-P_{n+1}(c)/P_n(c)$. Moreover,			
$$\widetilde{J}-cI=UL.$$				
\end{teo}
Observe that since $\un$ and $\widetilde\un$ are quasi-definite linear functionals, then so are $\widetilde \un^{(1)}$ and $\un^{(1)}$. Let $(\widetilde P^{(1)}_n)_{n\geq 0}$ and $(P_n^{(1)})_{n\geq 0}$ be their SMOP, respectively.
We are interested in analyzing the relation between $\widetilde \un^{(1)}$ and $\un^{(1)}$.

\begin{pro} The polynomial $\widetilde P^{(1)}_n(x)$ satisfies the following connection formula
$$(x-c)\widetilde P^{(1)}_{n-1}(x)=R_n(x)-\dfrac{P_{n+1}(c)}{P_n(c)}R_{n-1}(x),\quad n\geq 0,$$
where $R_n(x)=\dfrac{\un_0}{\widetilde \un_0}\left[(x-c)P^{(1)}_n(x)-P_{n+1}(x)\right].$	
\end{pro}
\begin{proof}
From \eqref{repChris}, we have
$$\prodint{\un_y,
	\frac{(x-c)\widetilde{P}_n(x)-(y-c)\widetilde{P}_n(y)}{x-y}}=\un_0P_n^{(1)}(x)-\un_0\dfrac{P_{n+1}(c)}{P_n(c)}P_{n-1}^{(1)}(x).$$
On the other hand, the left hand side of the previous equation is equivalent to
\begin{multline*}
\prodint{\un_y,\frac{(x-c)\widetilde{P}_n(x)-(y-c)\widetilde{P}_n(x)+(y-c)\widetilde{P}_n(x)-(y-c)\widetilde{P}_n(y)}{x-y}}=\\
\prodint{\un_y,\widetilde P_n(x)}+\prodint{\widetilde\un_y,\dfrac{\widetilde P_n(x)-\widetilde P_n(y)}{x-y}}=
\un_0\widetilde P_n(x)+\widetilde \un_0\widetilde P^{(1)}_{n-1}(x).
\end{multline*}
This yields
$$ \widetilde \un_0\widetilde P^{(1)}_{n-1}(x)=\un_0P_n^{(1)}(x)-\un_0\dfrac{P_{n+1}(c)}{P_n(c)}P_{n-1}^{(1)}(x)-\un_0\widetilde P_n(x).$$
Multiplying both hand sides of the above equation by $(x-c)$  and taking into account \eqref{repChris}, the statement follows.
\end{proof}	
\begin{pro}\label{pro5}
The polynomials $R_n(x)$ are co-recursive of parameter $\alpha:=- a_1\dfrac{\un_0}{\widetilde \un_0}$ with respect to the linear functional $\un^{(1)}$.
\end{pro}
\begin{proof}
From the three term recurrence relations (\ref{ttrr}) and (\ref{asociados}), we have
\begin{align*}
x\,P_{n}(x)&=P_{n+1}(x)+b_n\,P_{n}(x)+ a_{n}\,P_{n-1}(x),\quad n\geq 0,
\end{align*}
and
\begin{multline*}
x(x-c)\,P^{(1)}_{n-1}(x)=(x-c)P^{(1)}_{n}(x)+b_n\,(x-c)P^{(1)}_{n-1}(x)\\+ a_{n}\,(x-c)P^{(1)}_{n-2}(x),\quad n\geq 1,	
\end{multline*} thus
$$	x\,R_{n-1}(x)=R_{n}(x)+b_n\,R_{n-1}(x)+ a_{n}\,R_{n-2}(x),\quad n\geq0.$$
Taking into account that $(x-c)P_0^{(1)}(x)-P_1(x)=b_0-c$, and $(b_0-c)=\widetilde \un_0/\un_0\ne 0,$ then
$$R_1(x)=(x-b_1)R_0(x)-a_1R_{-1}(x)=P_1^{(1)}(x)+ a_1\dfrac{\un_0}{\widetilde \un_0},$$
and we get the result.
\end{proof}
Let $\un^\alpha$ be the quasi-definite linear functional such that $(R_n)_{n\geq 0}$ is its SMOP.  Since    $(R_n)_{n\geq 0}$ are co-recursive polynomials with respect to $\un^{(1)}$, then from \eqref{funccorre} we get
$$\un^\alpha=\dfrac{\un_0^\alpha}{\un^{(1)}_0}\left(\left(\un^{(1)}\right)^{-1}- a_1\dfrac{\un_0}{\un_0^{(1)}\widetilde \un_0}\deltan^\prime_0\right)^{-1}.$$
\begin{coro}\label{coro1}
The sequence $(\widetilde{P}^{(1)}_{n})_{n\geq 0}$ is orthogonal with respect to $(x-c)\un^\alpha$. In other words, $\widetilde\un^{(1)}$ is a canonical Christoffel transformation of~$\un^\alpha$.
\end{coro}
\begin{pro} Let $J$ and $J_\alpha$ be the monic Jacobi matrices associated with $(P_n)_{n\geq 0}$ and $(R_n)_{n\geq 0}$, respectively. If $J-cI$ has a LU factorization as in \eqref{LU}, then $J_\alpha-cI$ also has a LU factorization as follows
\begin{equation*}
J_\alpha-cI=:L_1U_1=\begin{pmatrix}
1           &            &&\\
\ell_{2}&1           &&\\
&\ell_{3}&1&\\
&            &\ddots&\ddots
\end{pmatrix}\begin{pmatrix}
\beta_{1}&1            &             &&\\
&\beta_{2}&1            &&\\
&             &\beta_{3}&\ddots &\\
&             &             &\ddots&\ddots
\end{pmatrix}.\end{equation*}	
Moreover, if $\widetilde J^{(1)}$ is the monic Jacobi matrix associated with $(\widetilde P_n^{(1)})_{n\geq 0}$, then $$\widetilde J^{(1)}-cI=U_1L_1.$$
\end{pro}
\begin{proof}
Notice that from the LU factorization of $J$
$$\begin{pmatrix}
1           &            &&\\
\ell_{2}&1           &&\\
&\ell_{3}&1&\\
&            &\ddots&\ddots
\end{pmatrix}\begin{pmatrix}
\beta_{1}&1            &             &&\\
&\beta_{2}&1            &&\\
&             &\beta_{3}&\ddots &\\
&             &             &\ddots&\ddots
\end{pmatrix}=\begin{pmatrix}
\beta_1 &1  &   &\\
a_2 &b_2-c& 1&\\
& a_3&b_3-c&\ddots\\
& & \ddots&\ddots
\end{pmatrix}.$$
On the other hand, from Proposition \ref{pro5} the monic Jacobi matrix associated with the SMOP $(R_n)_{n\geq 0}$ is given by
$$J_\alpha=\begin{pmatrix}
b_1+\alpha &1  &   &\\
a_2 &b_2& 1&\\
& a_2&b_3&\ddots\\
& & \ddots&\ddots
\end{pmatrix},\quad \text{with} \quad \alpha=- a_1\dfrac{\un_0}{\widetilde \un_0}.$$
Thus, the proof will be completed if we show that $b_1+\alpha=\beta_1+c.$ But
$$b_1+\alpha=b_1-\frac{ a_1}{(b_0-c)}=b_1+\frac{(c-b_1)P_1(c)-P_2(c)}{P_1(c)}=c+\beta_1.$$
From here, we get the result. The second part is a straightforward consequence of Theorem~\ref{teo1} and Corollary \ref{coro1}.
\end{proof}
\begin{rem}
Observe that LU factorization for $J_\alpha-cI$ is precisely the LU factorization for $J-cI$ but in each matrix, $L$ and $U$, we have removed its first column and its first row. In other words
$$L_1=\Lambda L\Lambda^\top,\quad U_1=\Lambda U\Lambda^\top,$$
where $\Lambda$ is the shift matrix given by
\begin{equation}\label{Lambda}
\Lambda=\begin{pmatrix}
0&1&0&&\\
0&0&1&0&\\
\vdots&\ddots&\ddots&\ddots&\ddots
\end{pmatrix}.
\end{equation}
\end{rem}
\subsection{Geronimus transformation}

Let $\vn$ be a quasi-definite linear functional and let $(P_n)_{n\geq 0}$ be its corresponding SMOP. If $c$ is a complex number, the linear functional $\widehat \vn$ defined by $(x-c)\widehat\vn=\vn$ is said to be the canonical Geronimus transformation of the linear  functional $\vn$.
Observe that $\widehat \vn$ is not uniquely defined since its first moment is arbitrary. The explicit expression of $\widehat \vn$ is given by \cite{Yoon02}
\begin{equation}\label{explicit}
\widehat \vn=(x-c)^{-1}\vn+\widehat\vn_0\deltan_c.
\end{equation}
Suppose that $\widehat\vn$ is also quasi-definite and let $(\widehat P_n)_{n\geq 0}$ be its SMOP.  It is well known that $(P_n)_{n\geq 0}$ and $(\widehat P_n)_{n\geq 0}$ are related by \cite{Bue1,Yoon02}

\begin{equation}\label{U}
\widehat{P}_n(x)=P_{n}(x)+\ell_nP_{n-1}(x), \quad n\geq 1,
\end{equation}
where
\begin{equation}\label{ln}
\ell_n=-\dfrac{\vn_0P^{(1)}_{n-1}(c)+\widehat\vn_0 P_n(c)}{\vn_0P^{(1)}_{n-2}(c)+\widehat\vn_0 P_{n-1}(c)},\quad n\geq 1.
\end{equation}
Thus a  necessary and sufficient  condition on $\widehat\vn$ to be a quasi-definite functional is
\begin{equation}\label{condition}
\widehat\vn_0\ne -\dfrac{\vn_0P^{(1)}_{n-1}(c)}{P_n(c)},\quad \text{for all}\  n\geq 1. 	\end{equation}

A second equation relating $(P_n)_{n\geq 0}$ and   $(\widehat P_n)_{n\geq 0}$ is the following one (see \cite{Yoon02})
\begin{equation}\label{rel}
(x-c)P_n(x)=\widehat P_{n+1}(x)+\beta_n\widehat P_n(x),\quad n\geq0,
\end{equation}
where $\beta_n=-\widehat P_{n+1}(c)/\widehat P_n(c).$ With this in mind we have the following relation between the corresponding monic Jacobi matrices.

\begin{teo}[\cite{Bue1,DereM,Dere,Yoon02}]
Let $J$ and $\widehat J$ be the monic Jacobi matrices associated with  $\vn$ and  $\widehat{\vn}$, respectively. If $\widehat\vn_0$ satisfies \eqref{condition}, then $J-cI$ has UL factorization. Indeed,
\begin{equation} \label{UL}
J-cI:=UL:=\begin{pmatrix}
\beta_{0}&1            &             &&\\
&\beta_{1}&1            &&\\
&             &\beta_{2}&\ddots &\\
&             &             &\ddots&\ddots
\end{pmatrix}\begin{pmatrix}
1           &            &&\\
\ell_{1}&1           &&\\
&\ell_{2}&1&\\
&            &\ddots&\ddots
\end{pmatrix},\end{equation}
where $L$ is a lower bidiagonal matrix with 1's as diagonal entries and $U$ is an upper bidiagonal matrix with $\beta_n=-\widehat P_{n+1}(c)/\widehat P_n(c)$. Moreover			
$$\widehat{J}-cI=LU.$$				
\end{teo}
Observe that the $UL$ factorization
depends on the choice  of $\widehat\vn_0$ since   $\beta_0=~\vn_0/\widehat\vn_0$.\\

Since we  assume that $\vn$ and $\widehat\vn$ are quasi-definite linear functionals, then    $\vn^{(1)}$ and $\widehat \vn^{(1)}$ are also quasi-definite linear functionals. Let   $(P_n^{(1)})_{n\geq 0}$ and $(\widehat P^{(1)}_n)_{n\geq 0}$ be their SMOP, respectively.
Now we are interested in analyzing the relation between $ \vn^{(1)}$ and $\widehat\vn^{(1)}$.

\begin{pro} The polynomial $\widehat P^{(1)}_n(x)$ satisfies the following connection formula
\begin{equation}\label{gero1}
(x-c)\widehat P^{(1)}_{n-1}(x)=S_n(x)+\ell_nS_{n-1}(x),\quad n\geq 1,
\end{equation}
where $S_n(x)= P_{n}(x)+\dfrac{\vn_0}{\widehat\vn_0}P^{(1)}_{n-1}(x).$	
\end{pro}
\begin{proof}
From \eqref{asociadosrepr} with $k=1$ and \eqref{explicit} we have
\begin{align*}
&\widehat\vn_0(x-c)\widehat{P}_{n-1}^{(1)}(x)\\
&=\prodint{\vn_y,\frac{(y-c)\widehat{P}_n(x)+(c-x)\widehat{P}_n(y)+(x-y)\widehat{P}_n(c) }{(x-y)(y-c)}} +\widehat\vn_0(\widehat{P}_n(x)-\widehat{P}_n(c))\\
&=\prodint{\vn_y,\frac{\widehat{P}_n(x)-\widehat{P}_n(y)}{(x-y)}}-\prodint{\vn_y,\frac{\widehat{P}_n(y)-\widehat{P}_n(c)}{(y-c)}}+\widehat\vn_0(\widehat{P}_n(x)-\widehat{P}_n(c))\\
&=\prodint{\vn_y,\frac{\widehat{P}_n(x)-\widehat{P}_n(y)}{(x-y)}}-\prodint{\widehat\vn_y,\widehat P_n(y)}+\widehat\vn_0\widehat{P}_n(x).\end{align*}
Taking into account \eqref{U} and \eqref{ln} we get the statement.
\end{proof}

\begin{pro}
The SMOP $(S_n)_{n\geq0}$ is the co-recursive SMOP of parameter $\alpha:=-\dfrac{\vn_0}{\widehat \vn_0}$ with respect to linear  functional $\vn$. Moreover, the linear functional $\vn^\alpha$ associated with the above sequence  can be written as

$$\vn^\alpha=\dfrac{\vn_0^\alpha}{\vn_0}\left(\vn^{-1}-\dfrac{1}{\widehat\vn_0}\deltan^\prime_0\right)^{-1}.$$		
\end{pro}
\begin{proof}
It is a direct consequence of the recurrence relation \eqref{ttrr1} and the fact that $S_1(x)=P_1(x)+\dfrac{\vn_0}{\widehat \vn_0}$. Moreover, the expression of $\vn^\alpha$ can be obtained from~\eqref{funccorre}.
\end{proof}

From \eqref{rel} we obtain the following connection formula
\begin{equation}\label{gero2}
S_n(x)=\widehat P^{(1)}_{n}(x)+\beta_n\widehat P^{(1)}_{n-1}(x),\quad n\in\N.
\end{equation}
As a consequence,  we can state the following proposition.

\begin{pro}\label{pro6} Let $J$ and $J_\alpha$ be the monic Jacobi matrices associated with $(P_n)_{n\geq 0}$ and $(S_n)_{n\geq0},$ respectively. Then
\begin{enumerate}
	\item $\widehat\vn^{(1)}$ is a Christoffel transformation of $\vn^\alpha,$ that is $\widehat\vn^{(1)}=(x-c)\vn^\alpha$.
	
	\item	If $J-cI$ has an UL factorization as in \eqref{UL}, then $J_\alpha-cI$ has the following LU factorization
	\begin{equation*}
	J_\alpha-cI=:\widehat L\widehat U=\begin{pmatrix}
	1           &            &&\\
	\beta_{1}&1           &&\\
	&\beta_{2}&1&\\
	&            &\ddots&\ddots
	\end{pmatrix}\begin{pmatrix}
	\ell_{1}&1            &             &&\\
	&\ell_{2}&1            &&\\
	&             &\ell_{3}&\ddots &\\
	&             &             &\ddots&\ddots
	\end{pmatrix}.\end{equation*}	
	
	Moreover,  if $\widehat J^{(1)}$ is the monic Jacobi matrix associated with $(\widehat P_n^{(1)})_{n\geq 0}$, then $$\widehat J^{(1)}-cI=\widehat U \widehat L.$$
\end{enumerate}
\end{pro}
\begin{rem}
Observe that
$$\Lambda L=\widehat U,\quad U\Lambda^\top=\widehat L,$$ where $\Lambda$ is the semi-infinite shift matrix given in \eqref{Lambda}.
\end{rem}

\section{Examples}\label{section4}
\begin{exa}
Let $(U_n)_{n\geq0}$ be the sequence of  monic Chebyshev polynomials of the second kind which are orthogonal with respect to the positive definite linear functional $\un$ defined by
\begin{equation}\label{cheby21}
\prodint{\un,p(x)}=\int_{-1}^1p(x)(1-x^2)^{1/2}\,dx,\quad p(x)\in\mathbb{P}.
\end{equation}
They satisfy the following properties:
\begin{enumerate}
	\item Recurrence relation.
	\begin{align*}
	xU_n(x)&=U_{n+1}(x)+\dfrac{1}{4}U_{n-1}(x),\quad n\geqslant0,\\
	U_{-1}(x)&=0,\quad U_0(x)=1.
	\end{align*}
	\item Values at $x=0$.
	$$U_n(0)=\begin{cases}
	1,&\textnormal{if $n=0$},\\
	0,& \textnormal{ if $n$ odd},\\
	\dfrac{(-1)^{\frac{n}{2}}}{2^{n}},&\textnormal{otherwise},
	\end{cases},\quad U_n'(0)=\begin{cases}
	0,&\textnormal{if $n$ even},\\
	\dfrac{(-1)^{\frac{n-1}{2}}}{2^{n}}(n+1),&\textnormal{if $n$ odd}.
	\end{cases}$$
	\item Normalization.
	$$\prodint{\un,U_n(x)U_m(x)}=2^{-(n+m+1)}\pi\delta_{n,m}.$$
\end{enumerate}
According to the definition of $(U_n^{(1)})_{n\geqslant0}$, it turns out that, for $n\geq0$, $U_n^{(1)}(x)=U_n(x)$. 	
From \eqref{deri}, \eqref{uminus} and \eqref{cheby21} we get that
\begin{equation}\label{cheb2uinv}
\prodint{\un^{-1},p(x)}=-\frac{1}{\pi^2}\int_{-1}^1\Biggl(\frac{p(x)-p(0)-p'(0)x}{x^2}\Biggr)(1-x^2)^{1/2}\, dx+\frac{2}{\pi}p(0).
\end{equation}
Corollary \ref{pnminusexpre} yields, for $n\geqq 2$,
$$d_n^*=\begin{cases}
\dfrac{-2^{1-2n}n}{\pi^2},&\textnormal{if $n$ even},\\[10pt]
-\dfrac{2^{1-2n}(n+1)}{\pi^2},&\textnormal{if $n$ odd},
\end{cases}$$
and
\begin{align*}
U^-_n(x)&=U_n(x)+\alpha_{n,n-2}U_{n-2}(x),\quad \alpha_{n,n-2}=\begin{cases}
\dfrac{n+2}{4n},&\textnormal{if $n$ even},\\[10pt]
\dfrac{1}{4},&\textnormal{if $n$ odd},
\end{cases}\quad n\geq2,\\
U^-_0(x)&=1,\quad U^-_1(x)=x.
\end{align*}
From Corollary \ref{ttrrpminus}
$$xU^-_n(x)=U^-_{n+1}(x)+ a_n^-U_{n-1}^-(x),\quad  a_n^-=\begin{cases}
-\dfrac{1}{4},&\textnormal{if $n=1$},\\[10pt]
\dfrac{n+2}{4n},&\textnormal{if $n$ even},\\[10pt]
\dfrac{n-1}{4(n+1)},&\textnormal{otherwise}.
\end{cases}$$
In terms of the Stieltjes function it is not difficult to check from \eqref{cheb2uinv} that
$$ \Su_{\un^{-1}}(z)=\frac{2}{\pi z}-\frac{1}{\pi^2 z^2} \Su_{\un}(z).$$ 	
\end{exa}
\begin{exa}
Let $(T_n)_{n\geq0}$ be the sequence of  monic Chebyshev polynomials of the first kind which are orthogonal with respect to the positive definite linear functional $\vn$ defined by
\begin{equation}\label{chebyt}
\prodint{\vn,p(x)}=\int_{-1}^1p(x)(1-x^2)^{-1/2}\,dx,\quad p(x)\in\mathbb{P}.
\end{equation}
They satisfy the following properties:
\begin{enumerate}
	\item Recurrence relation.
	\begin{align*}
	xT_n(x)&=T_{n+1}(x)+ a_{n}T_{n-1}(x),\quad n\geqslant0,\\
	T_{-1}(x)&=0,\quad T_0(x)=1,
	\end{align*}
	where $ a_{1}=1/2$ and $ a_{n}=1/4$ for $n\geqq2$.
	\item Values at $x=0$.
	$$T_n(0)=\begin{cases}
	1,&\textnormal{if $n=0$},\\
	0,& \textnormal{ if $n$ odd},\\
	(-1)^{\frac{n}{2}}2^{1-n},&\textnormal{otherwise},
	\end{cases},\quad T_n'(0)=\begin{cases}
	0,&\textnormal{if $n$ even},\\
	(-1)^{\frac{n-1}{2}}2^{1-n}n,&\textnormal{if $n$ odd}.
	\end{cases}$$
	\item Normalization.
	$$\prodint{\vn,T_n(x)T_m(x)}=\begin{cases}
	2^{1-2n}\pi\delta_{n,m},&\textnormal{if $n>0$},\\
	\pi,&\textnormal{if $m=n=0$}.
	\end{cases}$$
\end{enumerate}
According to the definition of $(T_n^{(1)})_{n\geqslant0}$, it turns out that, for $n\geq0$, $T_n^{(1)}(x)=U_n(x)$ and $\vn^{(1)}=\un$ given in \eqref{cheby21}. From \eqref{deri}, \eqref{uminus}, \eqref{cheby21} and \eqref{chebyt} we get that
$$\prodint{\vn^{-1},p(x)}=-\frac{1}{\pi^2}\int_{-1}^1\Biggl(\frac{p(x)-p(0)-p'(0)x}{x^2}\Biggr)(1-x^2)^{1/2}\, dx+\frac{1}{\pi}p(0).$$
Corollary \ref{pnminusexpre} yields, for $n\geqq 2$,
$$d_n^*=\begin{cases}
\frac{2^{3-2n}(1-n)}{\pi^2},&\textnormal{if $n$ even},\\[5pt]
-\frac{2^{3-2n}n}{\pi^2},&\textnormal{if $n$ odd,}
\end{cases}$$
and
\begin{align*}
T^-_n(x)&=U_n(x)+\alpha_{n,n-2}U_{n-2}(x),\quad \alpha_{n,n-2}=\begin{cases}
\frac{n+1}{4(n-1)},&\textnormal{if $n$ even},\\[5pt]
\frac{1}{4},&\textnormal{if $n$ odd},
\end{cases}\quad n\geq2,\\
T^-_0(x)&=1,\quad T^-_1(x)=x.
\end{align*}
Finally, according to Corollary \ref{ttrrpminus}
$$xT^-_n(x)=T^-_{n+1}(x)+ a_n^-T_{n-1}^-(x),\quad  a_n^-=\begin{cases}
-\frac{1}{2},&\textnormal{if $n=1$},\\[5pt	]
\frac{n+1}{4(n-1)},&\textnormal{if $n$ even},\\[5pt]
\frac{n-2}{4n},&\textnormal{otherwise}.
\end{cases}$$
Notice that
$T_n(\alpha_k)=T_n^{(1)}(\alpha_k)=T_n^-(\alpha_k)$ if $\alpha_k=\cos\theta_k$ with $\theta_k=\frac{k\pi}{n-1}$, $k=1,\ldots,n-2$, and
$$ \Su_{\vn^{-1}}(z)=\frac{1}{\pi z}-\frac{1}{\pi^2 z^2} \Su_{\vn^{(1)}}(z)$$
\end{exa}

\begin{exa}Let $(L_n^{\alpha+1}(x))_{n\geq0}$ be the sequence of monic Laguerre polynomials of parameter $\alpha+1$ with $\alpha>-1,$ which are orthogonal with respect to the positive definite linear functional $\vn$ defined by

$$\prodint{\vn,p(x)}=\int_{0}^\infty p(x)x^{\alpha+1}e^{-x}dx,\quad p(x)\in\mathbb{P}.$$	
They satisfy the following properties
\begin{enumerate}
	\item Recurrence relation.
	\begin{equation*}
	xL_n^{\alpha+1}(x)=L_{n+1}^{\alpha+1}(x)+(2n+\alpha+2)L_n^{\alpha+1}(x)+n(n+\alpha+1)L_{n-1}^{\alpha+1}(x), \quad n\geq 0,
	\end{equation*}
	with $L_{0}^{\alpha+1}(x)=1$ and  $L_{-1}^{\alpha+1}(x)=0.$
	\item $\left(\dfrac{d^i}{x^i}L_n^{\alpha+1}\right)(0)=(-1)^{n+i}\dfrac{n! \Gamma(\alpha+n+2)}{(n-i)! \Gamma(\alpha+i+2)}.$
	\item $\prodint{\vn,L_n^{\alpha+1}(x)L_m^{\alpha+1}(x)}=n! \Gamma(n+\alpha+2)\delta_{n,m}.$
\end{enumerate}

On the other hand, in \cite{AskWim1984} the authors studied the first kind Laguerre polynomials which are denoted by $(L^{\alpha+1}_{n}(x,1))_{n\geq0}$. In particular it was proved that these polynomials are orthogonal with respect to the positive definite functional $\vn^{(1)}$ defined by

$$\prodint{\vn^{(1)},p(x)}=\int_{0}^\infty p(x)\dfrac{x^{\alpha+1}e^{-x}}{\left|\Psi(1,-\alpha,xe^{-\pi i})\right|^2}dx,$$where$$\Psi(c,a,x)=\dfrac{1}{ \Gamma(c)}\int_0^{\infty e^{(3\pi/4)i}}t^{c-1}(1+t)^{a-c-1}e^{-xt}dt,$$
$\operatorname{Re}(c)>0,\ \ $ $-\pi/2<3\pi/4+\operatorname{arg}x<\pi/2.$

The monic associated polynomials of the first kind $(L^{\alpha+1}_{n}(x,1))_{n\geq0}$ satisfy the following properties
\begin{enumerate}
	\item Explicit formula.
	\begin{multline*}
	L_{n}^{\alpha+1}(x,1)=(-1)^n(n+1)(\alpha+3)_n\times\\\sum_{k=0}^n\dfrac{(-n)_kx^k}{(k+1)!(\alpha+3)_k}\times {}_{3}{F}_{2}\left(\begin{array}{c}{k-n,\ 1, \ \alpha+2} \\ {\alpha+k+3,\ k+2}\end{array} ; 1\right).
	\end{multline*}
	
	\item$L^{\alpha+1}_{n}(0,1)=\dfrac{(-1)^n}{\alpha+1}\left[(\alpha+2)_{n+1}-(n+1)!\right].$
	\item $\prodint{\vn^{(1)},L_n^{\alpha+1}(x,1)L_m^{\alpha+1}(x,1)}=(n+1)! \Gamma(n+\alpha+3)\delta_{n,m}.$
\end{enumerate}
It is straightforward to check that the linear functional $\vn^{-1}$ is  quasi-definite but not positive-definite. Moreover, from \eqref{deri} and \eqref{uminus}
\begin{multline*}
\prodint{\vn^{-1},p}=\dfrac{-(\alpha+1)}{ \Gamma(\alpha+2) \Gamma(\alpha+3)}\int_{0}^\infty \left(\dfrac{p(x)-p(0)-xp^\prime(0)}{x^2}\right)\dfrac{x^{\alpha+1}e^{-x}}{\left|\Psi(1,-\alpha,xe^{-\pi i})\right|^2}dx\\+\dfrac{p(0)}{ \Gamma(\alpha+2)}-\dfrac{p^\prime(0)}{ \Gamma(\alpha+1)}.
\end{multline*}
If $( \left(L_n^{\alpha+1}\right)^-(x))_{n\geq0}$ is the corresponding SMOP, using the previous properties and \eqref{inv}  we get the following connection formula 	for $n\geq 2$
$$\left(L_n^{\alpha+1}\right)^-(x)=L_{n}^{\alpha+1}(x,1)+2(n+\alpha+2)L_{n-1}^{\alpha+1}(x,1)+(n+\alpha+1)(n+\alpha+2)L_{n-2}^{\alpha+1}(x,1),$$
as well as $\left(L_1^{\alpha+1}\right)^-(x)=x+\alpha+2$ and  $\left(L_0^{\alpha+1}\right)^-(x)=1.$\\
With this in mind and \eqref{coe}, the polynomials $( \left(L_n^{\alpha+1}\right)^-(x))_{n\geq0}$ satisfy the following  three term recurrence relation
$$x\left(L_n^{\alpha+1}\right)^-(x)=\left(L_{n+1}^{\alpha+1}\right)^-(x)+b^-_{n}\left(L_{n}^{\alpha+1}\right)^-(x)+ a^-_{n}\left(L_{n-1}^{\alpha+1}\right)^-(x),$$
with
\begin{align*}
b^-_{n}&=(2n+\alpha+2),\quad n\geq 0,\\[10pt]
a^-_{n}&=(n-1)(n+\alpha+2), \quad  n\geq 2,\quad  a^-_1=-(\alpha+2)(\alpha+3).
\end{align*}

Next, let $\widehat\vn$ be the linear functional defined by the Geronimus transformation $x\widehat\vn=\vn$ with $\widehat\vn_0= \Gamma(\alpha+1)$.    Then from \eqref{explicit}
\begin{align*}
\prodint{\widehat\vn,p(x)}&=\int_{0}^\infty\left(p(x)-p(0)\right)x^\alpha e^{-x}\,dx+p(0)\int_0^\infty x^\alpha e^{-x}\,dx\\
&=\int_{0}^\infty p(x)x^{\alpha}e^{-x}dx,\quad p(x)\in\mathbb{P},
\end{align*}

where its corresponding sequence of monic orthogonal polynomials is   the sequence of monic  Laguerre polynomials of parameter $\alpha$. If $J_{\alpha+1}$ and $J_{\alpha}$ are the monic Jacobi matrices associated with $\vn$ and $\widehat\vn$ respectively,  then  $J_{\alpha+1}$ has UL factorization as in \eqref{UL} with $\beta_n=\alpha+n+1$, $n\geq 0$, $\ell_n=n$, $n\geq 1$, and  $ J_{\alpha}=LU$. From here we get the well-known formulas
\begin{align*}
L_n^{\alpha}(x)&=L_n^{\alpha+1}(x)+nL_{n-1}^{\alpha+1}(x),\\
xL_n^{\alpha+1}(x)&=L_{n+1}^{\alpha}(x)+(\alpha+n+1)L_{n}^{\alpha}(x).
\end{align*}
Let $(L_n^\alpha(x,1))_{n\geq0}$ be the SMOP associated with $\widehat\vn^{(1)}.$ Observe that there 	is no  relation between $\vn^{(1)}$ and $\widehat\vn^{(1)}$ of type $x\widehat\vn^{(1)}=\vn^{(1)}.$ However, from \eqref{gero1} and \eqref{gero2}  we get the following relations
\begin{align*}
xL_{n-1}^{\alpha}(x,1)&=L_n^{\alpha}(x)+(\alpha+1)L_{n-1}^{\alpha+1}(x,1)+n(\alpha+1)L_{n-2}^{\alpha+1}(x,1),\\
(\alpha+1)L_{n-1}^{\alpha+1}(x,1)&=L_{n}^{\alpha}(x,1)+(\alpha+n+1)L_{n-1}^{\alpha}(x,1)-L_n^{\alpha+1}(x).
\end{align*}

Let $ S_n(x):=L^{\alpha+1}_{n}(x)+(\alpha+1)L^{\alpha+1}_{n-1}(x,1)$ be the co-recursive polynomials of parameter $-(\alpha+1)$ which are orthogonal with respect to  (see \cite{Letessier}),

$$\prodint{\un,p(x)}=\int_{0}^\infty p(x)\frac{x^{\alpha+1} \mathrm{e}^{-x}}{\left|1-(\alpha+1) \Psi\left(1,-\alpha ; x \mathrm{e}^{-\mathrm{i} \pi}\right)\right|^{2}}dx.$$
From Proposition \ref{pro6}, we get that $\widehat\vn^{(1)}=x\un$. Thus
\begin{align*}
\prodint{\widehat\vn^{(1)},p(x)}&=\int_{0}^\infty p(x)\frac{x^{\alpha+2} \mathrm{e}^{-x}}{\left|1-(\alpha+1) \Psi\left(1,-\alpha ; x \mathrm{e}^{-\mathrm{i} \pi}\right)\right|^{2}}dx\\&=\int_{0}^\infty p(x)\frac{x^{\alpha} \mathrm{e}^{-x}}{\left| \Psi\left(1,1-\alpha ; x \mathrm{e}^{-\mathrm{i} \pi}\right)\right|^{2}}dx,
\end{align*}

that is the
expected result. The last equality on the right hand side is a straightforward consequence of the following identity (see \cite[Sect. 6.6 eq. 7]{Erd1953} )
$$(a-c)\Psi(c,a,x)-x\Psi(c,a+1,x)+\Psi(c-1,a,x)=0.$$

\end{exa}

\section{Conclusions and further remarks}
In this contribution we have analyzed associated polynomials of the first kind with respect to linear functionals defined as linear spectral transformations of a given linear functional. We have focused our attention on canonical Christoffel and Geronimus transformations which constitute a generating system of the set of linear spectral transformations. Explicit expressions for the sequences of orthogonal polynomials for such a kind of transformations as well as the relation between the corresponding tridiagonal matrices are obtained.


\section{Acknowledgements}
The work of the second author (FM) has been supported by Agencia Estatal de Investigaci\'on of Spain, grant PGC2018-096504-B-C33.

\end{document}